\documentclass[11pt,reqno]{amsart}
\usepackage{amssymb,mathrsfs,graphicx}
\usepackage{color}

\catcode`\@=11 \@addtoreset{equation}{section}

\catcode`\@=12

\newtheorem{Theorem}{Theorem}[section]
\newtheorem{Proposition}{Proposition}[section]
\newtheorem{Lemma}{Lemma}[section]
\newtheorem{Corollary}{Corollary}[section]
\newtheorem{Remark}{Remark}[section]

\newcommand{\bTheorem}[1]{
\begin{Theorem} \label{T#1} }
\newcommand{\eT}{\end{Theorem}}

\newcommand{\bProposition}[1]{
\begin{Proposition} \label{P#1}}
\newcommand{\eP}{\end{Proposition}}

\newcommand{\bLemma}[1]{
\begin{Lemma} \label{L#1} }
\newcommand{\eL}{\end{Lemma}}

\newcommand{\bCorollary}[1]{
\begin{Corollary} \label{C#1} }
\newcommand{\eC}{\end{Corollary}}

\newcommand{\bRemark}[1]{
\begin{Remark} \label{R#1} }
\newcommand{\eR}{\end{Remark}}

\newcommand{\bFormula}[1]{
\begin{equation} \label{#1}}
\newcommand{\eF}{\end{equation}}

\newcommand{\Ov}[1]{\overline{#1}}

\newcommand{\DC}{C^\infty_c}

\newcommand{\vr}{\varrho}

\newcommand{\vu}{\vc{u}}
\newcommand{\vc}[1]{{\bf #1}}
\newcommand{\Div}{{\rm div}_x}

\newcommand{\Grad}{\nabla_x}

\newcommand{\tn}[1]{\mbox {\F #1}}
\newcommand{\dx}{{\rm d} {x}}
\newcommand{\dt}{{\rm d} t }

\newcommand{\intO}[1]{\int_{\Omega} #1 \ \dx}

\newcommand{\D}{{\mathcal D}}

\font\F=msbm10 scaled 1000
\newcommand{\R}{\mbox{\F R}}

\newcommand{\Del}{\Delta_x}

\definecolor{Cgrey}{rgb}{0.85,0.85,0.85}
\definecolor{Cblue}{rgb}{0.50,0.85,0.85}
\definecolor{Cred}{rgb}{1,0,0}
\definecolor{fancy}{rgb}{0.10,0.85,0.10}

\newcommand\Cbox[2]{%
    \newbox\contentbox%
    \newbox\bkgdbox%
    \setbox\contentbox\hbox to \hsize{%
        \vtop{
            \kern\columnsep
            \hbox to \hsize{%
                \kern\columnsep%
                \advance\hsize by -2\columnsep%
                \setlength{\textwidth}{\hsize}%
                \vbox{
                    \parskip=\baselineskip
                    \parindent=0bp
                    #2
                }%
                \kern\columnsep%
            }%
            \kern\columnsep%
        }%
    }%
    \setbox\bkgdbox\vbox{
        \color{#1}
        \hrule width  \wd\contentbox %
               height \ht\contentbox %
               depth  \dp\contentbox
        \color{black}
    }%
    \wd\bkgdbox=0bp%
    \vbox{\hbox to \hsize{\box\bkgdbox\box\contentbox}}%
    \vskip\baselineskip%
}

\date{}

\begin{document}


\title[On weak solutions to Euler systems with non-local interactions]{Weak solutions for Euler systems with non-local interactions}

\author[Carrillo]{Jos\'{e} A. Carrillo}
\address[Jos\'{e} A. Carrillo]{\newline Department of Mathematics
    \newline Imperial College London, London SW7 2AZ, United Kingdom}
\email{carrillo@imperial.ac.uk}

\author[Feireisl]{Eduard Feireisl }
\address[Eduard Feireisl]{\newline Institute of Mathematics of the Academy of Sciences of the Czech Republic\newline \v Zitn\' a 25, 115 67 Praha 1, Czech Republic}
\email{feireisl@math.cas.cz}

\author[Gwiazda]{Piotr Gwiazda}
\address[Piotr Gwiazda]{\newline Institute of Applied Mathematics and Mechanics, University of Warsaw \newline Banacha 2, 02-097 Warszawa, Poland}
\email{pgwiazda@mimuw.edu.pl}

\author[\' Swierczewska-Gwiazda]{Agnieszka \' Swierczewska-Gwiazda}
\address[Agnieszka \' Swierczewska-Gwiazda]{\newline Institute of Applied Mathematics and Mechanics, University of Warsaw \newline Banacha 2, 02-097 Warszawa, Poland}
\email{aswiercz@mimuw.edu.pl}

\maketitle

\begin{abstract}

We consider several modifications of the Euler system of fluid dynamics including its pressureless
variant driven by non-local interaction repulsive-attractive and alignment forces in the space dimension $N=2,3$. These models arise in the study of self-organisation in collective behavior modeling of animals and crowds.
We adapt the method of convex integration to show the existence of infinitely
many global-in-time weak solutions for any bounded initial data. Then we consider the class of \emph{dissipative} solutions satisfying, in addition, the associated global energy balance (inequality). We identify a large set of initial data for which the problem admits infinitely many dissipative weak solutions. Finally, we establish a weak-strong uniqueness principle for
the pressure driven Euler system with non-local interaction terms as well as for the  pressureless system with Newtonian interaction.

\end{abstract}

\bigskip


\section{Introduction}
\label{i}

We study various models arising in mathematical biology based on fluid dynamics approximations of Euler system type.
As singularities are expected to appear in a finite time no matter how smooth the initial data are, we consider the problems
in the framework of weak solutions. Adapting the method of convex integration we show existence of infinitely many weak
solutions for \emph{any} sufficiently smooth initial data.
The heart of the method is the so-called
oscillatory lemma shown in the seminal work by DeLellis and Sz\' ekelyhidi \cite{DelSze3}, adapted to the compressible Euler system by Chiodarolli \cite{Chiod} and to more general ``variable coefficients'' problems in \cite{ChiFeiKre, DoFeMa, EF2015_1, FGSG}.
On the other hand we show that the class of \emph{strong} solutions is robust in terms of a weak-strong uniqueness principle: weak and strong solutions emanating from the same initial data coincide as long as the latter exists.

The most general form of the class of problems we are interested in
reads, see \cite{CCTT}:
\begin{align}
\partial_t \vr(t,x) + \Div \Big( \vr(t,x) \vu(t,x) \Big) =& \, 0,
\label{i1}\\
\partial_t \Big( \vr(t,x) \vu(t,x) \Big) + \Div \Big( \vr(t,x) \vu(t,x)  \otimes \vu(t,x) \Big) = &- \Grad p(\vr(t,x)) + \left( 1 - H\left( |\vu(t,x) |^2 \right) \right) \vr(t,x) \vu(t,x)
\nonumber\\
&- \vr(t,x) \int_{\Omega} \Grad K(x - y) \vr(t,y) \ {\rm d}y
\label{i2}\\
&+
\vr(t,x) \int_{\Omega} \psi (x-y) \Big( \vu(t,y) - \vu(t,x) \Big) \vr(t,y) \ {\rm d}y ,
\nonumber
\end{align}
where $\vr = \vr(t,x)$ is the density, $\vu = \vu(t,x)$ the velocity field, $p = p(\vr)$ the pressure, and the $K$ and $\psi$ represent non-local interaction forces acting on the medium. More precisely, the kernel $K$ includes the repulsive-attractive interaction force between individuals \cite{CDMBC,CDP} and $\psi$ gives the local averaging measuring the consensus in orientation among individuals \cite{HT}. Note that (\ref{i1}), (\ref{i2}) reduces to the standard barotropic Euler system provided
$K = \psi = 0$, $H \equiv 1$.

These models are derived as hydrodynamic solutions to kinetic equations for swarming \cite{CDP,CKMT}. The fundamental \emph{Individual Based Models} (particle models) have been classically proposed by the applied mathematics and theoretical biology communities to describe the complicated patterns and coherent structures observed in animal behavior, see  \cite{HW,LR,CDF,Bir} and the references therein. Some minimal models for swarming include the 3-zone type models for which the main interactions taken into account between individuals are: An inner-core repulsion for comfort or collision avoidance, a global tendency to attract each other in such a way that compact groups are formed, and finally a general mechanism of gregariousness in which their orientation or velocity are influenced by their neighbors. Basic models can be found in the applied mathematics literature, see \cite{DCBC,CS1,CS2}, improved by including several other effects such as sensitivity regions, roosting forces, among others \cite{CKMT,CCR2,AIR,CCH} and the references therein. Repulsion-attraction interactions are assumed to be given by social forces imposed by the presence of an individual in the sensitivity region, in the simplest model introduced in \cite{DCBC} they correspond to a force field which is the gradient of a given radial potential $K(x)$ depending on the interparticle distance. Alignment interactions are modelled \cite{CS1,CS2,MT} by locally averaging the relative velocity with respect to neighbors with weights depending on the inter-particle distances and being normalized or not.

Apart from the aforementioned mechanisms, many models assume that particles are already moving with some cruise speed to avoid uninteresting steady situations. This effect is often imposed dynamically, that is, by imposing a term of friction between the particles whose unique speed root determines the asymptotic cruise speed of individuals \cite{DCBC}. All these effects being put together, the basic particle models look like
\begin{flalign}\label{micro}
\begin{split}
   \frac{dx_i}{dt} &= v_i \,, \\
   \frac{dv_i}{dt} &= v_i  - \alpha v_i \vert v_i\vert^2
   - \sum_{j \neq i} \nabla_{x} K(x_i-x_j) + \sum_{j } \psi(x_i-x_j) (v_j-v_i)\,,
\end{split}
\end{flalign}
where $i\in\{1,\dots,n\}$ and $\alpha>0$. Observe the asymptotic cruise speed is fixed to $\tfrac1{\sqrt{\alpha}}$. The passage $n\to\infty$ has been studied in the mean-field limit scaling by several authors, see \cite{CDP,HL,CCR,CCH} and the references therein. In the mean-field limit approximation the potential and the communication rate are scaled with the number of particles, and look at the limit of the empirical measure associated to the dynamical system \eqref{micro}, that is,
$$
\mu_n(t)= \frac1n \sum_{i=1}^n \delta_{(x_i(t),v_i(t))}
$$
defined as a probability measure in phase space $(x,v)\in \R^N\times\R^N$. Under certain assumptions on $K$ and $\psi$, it is proven that $\mu_n(t)$ converges as $n\to\infty$ to a solution of a Vlasov-like kinetic equation of the form
\begin{equation}
\partial_{t}f+v\cdot \nabla_{x}f +\mbox{div}_v\left(F[f] f +\left(1 - \alpha |v|^{2}\right)v f\right)
=0\,, \label{VP}
\end{equation}
with $F[f]$ being the nonlocal force field given by
$$
F[f] = -\nabla K\ast\vr + \int_{\R^N} \int_{\R^N} \psi(x-y) (w-v) f(t,x,w) dw \, dy
$$
and
\begin{equation*}
\vr(t,x):=\int_{\R^N} f(t,x,v) d v \,.
\end{equation*}
Here, $f(t,x,v): \R \times \R^N \times \R^N \rightarrow \R$ is the phase-space density, and $\vr(t,x)$ is the averaged (macroscopic) density.

In order to reduce complexity of the model, kinetic equations are usually being reduced to macroscopic equations by their closure or asymptotic limits. However, in the absence of collisions in the kinetic model, this could
represent a rather delicate task, and therefore special solutions of a particular form are sought. By imposing the so-called hydrodynamic or mono-kinetic ansatz, i.e., looking for distributional solutions to \eqref{VP} of the form
$$
f_h(t,x,v) = \vr(t,x) \delta(u-\vu(t,x))\,,
$$
one deduces that the pair $(\vr(t,x),\vu(t,x))$ must satisfy the hydrodynamic equations \eqref{i1}, \eqref{i2} with $H(s)=\alpha s^2$ without the pressure term. On the other hand, pressure terms can be easily derived from it if other closure assumptions are imposed dealing with particle models with noise or just by localizing the repulsive force terms via spatial scaling assumptions. The hydrodynamic equations \eqref{i1}, \eqref{i2} have recently been studied in one space dimension \cite{CCTT}, where interesting dichotomy results are available. If the attraction forces are not too strong, there are simple conditions separating two kinds of behavior: global existence of classical solutions and the blow-up of classical solutions in a finite time. The associated threshold conditions are known to appear also in other fluid mechanics problems of Euler type, see \cite{CCTT} and the references therein. Explicit global in time sufficiently regular solutions for this system can also be found in the form of flock solutions, see \cite{CHM,CDM}, the stability of which within the class of classical (strong) solutions represents an extremely interesting open problem.

In this paper, we show in a certain sense  negative results concerning stability of particular solutions. It turns out that the solutions must be sought in a stronger class than that of weak and/or dissipative solutions. We essentially show that there  are infinitely  many weak solutions for any initial data and that there is a vast class of velocity fields that gives rise to infinitely many admissible (dissipative) weak solutions. We may therefore infer that the class of weak solutions is not convenient for analysing certain qualitative properties such as stability and formation of the flock patterns. However, we also show that the strong solutions are robust in a larger class of all admissible (dissipative) weak solutions leading to the possibility of establishing certain stability results of flock solutions. 

Another naturally arising question would be whether the strong solutions were robust in even larger class of solutions, such as admissible measure valued solutions, see \cite{BDS, DiPeMa, GSW, DeStTz2012, Gw2} for studies in similar directions for different hyperbolic systems. 

For the sake of simplicity, we restrict ourselves to the spatially periodic boundary conditions, meaning $x \in \Omega$, where
\bFormula{i2a}
\Omega = \left( [-1, 1]|_{\{ -1,1 \}} \right)^N, \ N=2,3,
\eF
is the flat torus, noticing that the method can be adapted to general (smooth) bounded domains $\Omega \subset R^N$ endowed with the no-flux boundary conditions $\vu \cdot \vc{n}|_{\partial \Omega} = 0$. From now on, all convolutions must be understood in the flat torus.

The paper is organized as follows. In Section \ref{m}, we summarize our hypotheses imposed on the constitutive equations and state the main results of this work. In Section \ref{r}, the existence of weak solutions is rewritten in the form that allows to apply the abstract existence theorem proved in \cite{EF2015_1}. In Section \ref{d}, we introduce the dissipative solutions and show existence in this class for a large set of suitable initial data. Finally, we establish the weak-strong uniqueness principle in Sections \ref{WSU} and \ref{Poi}.


\section{Preliminaries \& Main results}
\label{m}

In this section, we list the principal hypotheses and state our main result concerning the existence of infinitely many weak solutions emanating from identical initial data.

\subsection{Hypotheses}

We suppose that the pressure $p$ is a continuously differentiable function of the density, specifically,
\bFormula{m1}
p \in C[0, \infty) \cap C^2(0, \infty), \ p(0) = 0.
\eF

\bRemark{m1}
There will be more restrictions imposed on the pressure required in the proof of the weak-strong uniqueness principle in Section \ref{WSU}. However, the mere hypothesis (\ref{m1}) is sufficient for the existence of infinitely many solutions stated below. In particular, the \emph{pressureless} Euler system with $p \equiv 0$ is included.
\eR

The kernels $K$ and $\psi$ represent the non-local interaction forces between individuals: $\nabla K$ encodes the repulsion-attraction forces while $\psi$ is the communication function between individuals for the local averaging of the alignment direction. We assume that the interaction potential and the communication rate are smooth enough such that the interaction force $\nabla K \ast \rho$ and the alignment force are regular $C^1(\Omega)$. This poses the following restrictions {
\bFormula{m2}
K \in C^2(\Omega), \ \psi \in C^1(\Omega), \ \psi \geq 0,\ \mbox{$K$ and $\psi$ are symmetric.}
\eF}
Note again that the situation $K \equiv 0$ and/or $\psi \equiv 0$ is allowed in the existence result stated below. In many real world applications, $K$ is assumed to be radial, repulsive in the short range and attractive in the long range, although these assumptions will not be needed in the existence theorem. The interaction potential and the communication rate are typically symmetric, $K(z)=K(-z)$ and $\psi(x)=\psi(-x)$, in the simplest particle models \cite{CS1,CFRT}. However, this symmetry assumption is often broken if orientational alignment mechanisms are taken into account such as vision or sensitivity regions for the individuals, see for instance \cite{MT,CFTV}. The symmetry of $\psi$ allows for simpler formulation of the dissipative nature of the alignment consensus for the velocity field as we will see below.

Finally, the function $H$ is related to ``friction'' similar to the Savage-Hutter model studied in \cite{FGSG}. Here, we assume{
\bFormula{m3}
H \in C^2([0, \infty)), \ 0 \leq H(Z) \leq H_\infty:=\lim_{Z\to\infty} H(Z), \ H'(Z) \geq 0\ \mbox{for all}\ Z \geq 0.
\eF}
The role of this term in the model is to fix an asymptotic speed for the individuals, i.e. a cruise speed. Then it is quite reasonable to assume that $H$ is an increasing function such that $1$ belongs to its range and as a consequence the equation $H(s^2)=1$ has a unique positive solution given by the cruise speed.

Even if some of the results below are true for more general assumptions on $K$, $\psi$ and $H$ than \eqref{m2}-\eqref{m3}, we prefer to state them in this way since they are the typical assumptions in their original application. We will point out after each theorem if some hypotheses are not needed or can be generalized with respect to \eqref{m2}-\eqref{m3}.

\subsection{Main results}

We consider the initial value problem for system \eqref{i1}, \eqref{i2}, \eqref{i2a} with the data
\bFormula{m4}
\vr(0, \cdot) = \vr_0, \ \vu(0, \cdot) = \vu_0 \ \mbox{in}\ \Omega.
\eF

We claim the following result regarding the existence of weak solutions.

\bTheorem{m1}
Let $T > 0$ be given and let $N=2,3$. Suppose that $p$, $\psi$, $K$ and $H$ satisfy the hypotheses \eqref{m1}--\eqref{m3}. Let the initial data
$\vr_0$, $\vu_0$ satisfy
$
\vr_0 \in C^2(\Omega), \ \vr_0 \geq \underline{\vr} > 0 \ \mbox{in}\ \Omega, \ \vu_0 \in C^3(\Omega; R^N).
$
Then the initial value problem \eqref{i1}, \eqref{i2}, \eqref{i2a} with initial data \eqref{m4} and periodic boundary conditions admits infinitely many weak solutions in the space-time cylinder $(0,T) \times \Omega$ belonging to the class
\bFormula{m5}
\vr \in C^2([0,T] \times \Omega), \ \vr > 0, \ \vu \in C_{\rm weak}([0,T]; L^2(\Omega; R^N)) \cap L^\infty((0,T) \times \Omega; R^N).
\eF
\eT

\bRemark{R2}\
\begin{itemize}
\item Here and hereafter,  \emph{weak solution} means a solution of \eqref{i1}, \eqref{i2} in the sense of generalized derivatives (distributions) in $[0,T) \times \Omega$ and $\vr\ge 0$, $\vr\in L^\infty(0,T;L^1(\Omega))$, 
$\vr|\vu|^2\in L^\infty(0,T;L^1(\Omega))$ and $P(\vr)\in L^\infty(0,T;L^1(\Omega))$. 

\item Hypothesis \eqref{m3} can be relaxed to
\[
0\leq H(Z) \leq c \left( 1 + |Z|^\alpha \right) \ \mbox{for some}\ \alpha < 1.
\]
This restriction is required for the set of the so-called subsolutions to be bounded, see Section \ref{r}.
\item The symmetry of $K$ and $\psi$ does not play a role in this result.
\end{itemize}
\eR

The proof of Theorem \ref{Tm1} based on the theory of convex integration will be given in the next section. We can still construct infinitely many solutions if we further restrict the concept of weak solutions by imposing that the dissipation of the total free energy of the system holds. More precisely, weak solutions satisfying
\begin{align*}
\intO{ \left[ \frac{1}{2} \vr |\vu|^2 + P(\vr) + \frac{1}{2} \vr K* \vr  \right] (\tau, \cdot) }
\leq &
\intO{ \left[ \frac{1}{2} \vr_0 |\vu_0|^2 + P(\vr_0) + \frac{1}{2} \vr_0 K* \vr_0  \right] }
\\
&+ \int_0^\tau \intO{ \vr |\vu|^2 \left(1 - H\left(|\vu|^2 \right)\right) } \ \dt
\\
&-\int_0^\tau\intO{ \int_\Omega \psi (x - y) \vr(t,x) \vr(t,y) \left| \vu(t,y) - \vu(t,x) \right|^2 \ {\rm d}y}\ \dt\,,
\end{align*}
for all $\tau>0$, where
\[
P(\vr) = \vr \int_1^\vr \frac{p(z)}{z^2} \ {\rm d}z
\]
is the pressure potential, are called dissipative weak solutions.
This concept of solution and the following result will be analysed in Section 4.

\bTheorem{d1}
Under the assumptions of Theorem \ref{Tm1},  given $T > 0$ and
$\vr_0 \in C^2(\Omega), \ \vr_0 > 0,$
there exists $\vu_0 \in L^\infty(\Omega; R^3)$ such that the initial value problem \eqref{i1}, \eqref{i2} and \eqref{m4}
admits infinitely many dissipative weak solutions in the space-time cylinder $(0,T) \times \Omega$.
\eT

\bRemark{R22}The symmetry of $\psi$ is not essential but leads to clearer dissipation properties. Hypothesis \eqref{m3} can be relaxed to the same property as in Remark \ref{RR2}.
\eR

We will finalize this work by showing weak-strong uniqueness results. Next result is shown in Section 5.

\bTheorem{e1}
Let $T > 0$ be given and let $N=2,3$. Suppose that  $p$, $\psi$, $K$ and $H$ satisfy the hypotheses \eqref{m1}--\eqref{m3} and in addition suppose that
\begin{equation} \label{tlak}
p(0) = 0, \
p'(\vr) > 0 \ \mbox{whenever}\ \vr > 0,\
\liminf_{\vr \to \infty} p'(\vr) > 0, \ \liminf_{\vr \to \infty} \frac{P(\vr)}{p(\vr)} > 0.
\end{equation}
Let $[\varrho,\vu]$ be a dissipative weak solution to system \eqref{i1}--\eqref{i2}, \eqref{m4} in $(0,T) \times \Omega$ 
 and let  $[r, \vc{U}]$  $r > 0$, be a Lipschitz (strong) solution 
%
%
with
\[
\varrho_0 = r(0, \cdot),\
\vu_0 = \vc{U}(0, \cdot).
\]
Then
\[
\varrho = r,\ \vu = \vc{U} \ \mbox{a.e. in}\ (0,T) \times \Omega.
\]
\eT

\begin{Remark} \label{rf1}\
\begin{itemize}
\item  We first point out that dissipative weak solutions in the previous theorem are distributional solutions not necessarily belonging to the class \eqref{m5}. In particular, the weak dissipative solutions for this result do not need to be bounded or bounded away from zero.

\item  Note that the condition \eqref{tlak} is satisfied for the ``isentropic'' pressure
\[
p(\vr) = a \vr^\gamma, \ \gamma \geq 1.
\]
In particular, the isothermal model derived by Karper, Mellet, and Trivisa \cite{KaMeTr} corresponding to $\gamma = 1$ is included.

\item  Hypothesis \eqref{m3} can be reduced to $H\in C^1([0,\infty))$ bounded, nonnegative and
\begin{equation} \label{HH1}
H'(Z)  \geq 0  \ \mbox{for all}\ Z \geq Z_0\,.
\end{equation}
\end{itemize}
\end{Remark}

We will be able to prove the same result in the case of the pressureless system with Newtonian interaction, see Theorem \ref{Te2} in Section 6. Finally, let us point out that a number of explicit flock solutions of the system \eqref{i1}-\eqref{i2}, \eqref{m4} without pressure in the whole space $(0,T) \times \R^N$ have been found, see \cite{CHM,CDM} and the references therein. They are given by density profiles $\bar\rho$ such that verify $\nabla K \ast \bar\rho=0$ in the support of $\bar\rho$ with the solutions given by 
$$
\rho(t,x)=\bar\rho(x-u_0 t) \qquad \mbox{and} \qquad u(t,x)=u_0\,,
$$
with $u_0$ arbitrary vector with unit norm. For instance, given $K(x)=\tfrac{|x|^2}2 - \frac{|x|^{2s-N}}{2s-N}$ with $0<s<1$, one can find that $\bar\rho$ is a Holder continuous compactly supported function that is $C^\infty$ inside its support. In the case of repulsive Newtonian, i.e. $s=1$, a flock solution is given by the characteristic function of a ball with a suitably chosen radius. By taking $K(x)=\tfrac{|x|^a}a - \frac{|x|^{2-N}}{2-N}$ with $2-N<a<2$ or $a>2$, there is a radially symmetric flock solution which is compactly supported and infinitely smooth inside its support with a jump discontinuity at the boundary of its support, see \cite{FH,FHK}. Also, in one dimension taking as potential $K(x)=\tfrac{x^2}2 - \log |x|$ the flock profile $\bar\rho$ is given by the semicircle law, see \cite{CFP}. Finding suitable classes of solutions to understand the stability of these flock solutions is a challenging open problem.


\section{Reduction to an abstract Euler system}
\label{r}

Our goal is to apply the abstract result \cite[Theorem 2.1]{EF2015_1} based on the oscillatory lemma of DeLellis and Sz\' ekelyhidi
\cite{DelSze3} for the \emph{incompressible} Euler system.
First, we introduce the notation
\[
\vc{v} \otimes \vc{w} \in R^{N \times N}_{\rm sym},\ [\vc{v} \otimes \vc{w}]_{i,j} = v_i v_j, \ \mbox{and}\
\vc{v}\odot \vc{w} \in R^{N \times N}_{{\rm sym},0}, \ \vc{v} \odot \vc{w} = \vc{v} \otimes \vc{w} - \frac{1}{N} \vc{v} \cdot \vc{w} \tn{I}.
\]
Following \cite{EF2015_1} we consider an abstract ``Euler problem'' in the form:

\medskip

{\it
Find a vector field $\vc{v} \in C_{\rm weak}([0,T]; L^2(\Omega; R^N))$ satisfying
\bFormula{E1}
\partial_t \vc{v} + \Div \left( \frac{ (\vc{v} + \vc{h}[\vc{v}] ) \odot (\vc{v} + \vc{h}[\vc{v}] ) }{r[\vc{v}]} + \tn{M}[\vc{v}] \right) = 0,\
\Div \vc{v} = 0
\ \mbox{in} \ \D'((0,T) \times \Omega; R^N),
\eF
\bFormula{E2}
\frac{1}{2} \frac{ | \vc{v} + \vc{h} [\vc{v}] |^2 }{r[\vc{v}]}  (t, x) = e[\vc{v}] (t,x) \ \mbox{for a.a.}\ (t,x) \in (0,T) \times \Omega,
\eF
\bFormula{E3}
\vc{v}(0, \cdot) = \vc{v}_0, \ \vc{v}(T, \cdot) = \vc{v}_T,
\eF
where $\vc{h}[\vc{v}]$, $r[\vc{v}]$, $\tn{H}[\vc{v}]$, and $e[\vc{v}]$ are given (nonlinear) operators.}

\medskip

Let $Q \subset (0,T) \times \Omega$ be an open set such that
\[
|Q| = |(0,T) \times \Omega|.
\]
The quantities $h[\vc{v}]$, $r = r[\vc{v}]$, $\tn{M} = \tn{M}[\vc{v}]$, $e = e[\vc{v}]$ appearing in (\ref{E1}--\ref{E3}) are assumed to be
non-local operators of casual type. More specifically, we say that
an operator
\[
b: C_{\rm weak}([0,T]; L^2(\Omega; R^N)) \cap L^\infty((0,T) \times \Omega; R^N) \to C_b (Q, R^M)
\]
is \emph{$Q-$continuous} if:

\begin{itemize}

\item $b$ maps bounded sets in $L^\infty((0,T) \times \Omega; R^N)$ on bounded sets in $C_b (Q, R^M)$;

\item $b$ is continuous, specifically,

\bFormula{E4a}
\begin{array}{c}
b[\vc{v}_n] \to b[\vc{v}] \ \mbox{in} \ C_b(Q; R^M) \ \mbox{(uniformly for $(t,x) \in Q$ ) } \\ \\
\mbox{whenever}\\ \\
\vc{v}_n \to \vc{v} \ \mbox{in}\ C_{\rm weak}([0,T]; L^2(\Omega; R^N)) \ \mbox{and weakly-(*) in} \ L^\infty((0,T) \times \Omega; R^N);
\end{array}
\eF

\item $b$ is causal (non-anticipative), meaning
$$
\vc{v}(t, \cdot) = \vc{w}(t,\cdot) \ \mbox{for}\ 0 \leq t \leq \tau \leq T
\ \mbox{implies} \ b[\vc{v}] = b[\vc{w}] \ \mbox{in}\  \left[ (0, \tau] \times \Omega \right] \cap Q.
$$

\end{itemize}

Finally, following DeLellis and Sz\' ekelyhidi \cite{DelSze3}, we introduce the space $X_0$ of \emph{subsolutions} associated to problem (\ref{E1}--\ref{E3}):
\bFormula{E6}
X_0 = \left\{ \vc{v} \ \Big| \ {\vc{v}} \in C_{\rm weak}([0,T]; L^2(\Omega; R^N)) \cap L^\infty((0,T) \times \Omega;R^N) \vphantom{\frac{1}{2}},\
\vc{v}(0, \cdot) = \vc{v}_0 , \ \vc{v}(T, \cdot) = \vc{v}_T, \right.
\eF
\[
\partial_t \vc{v} + \Div \tn{F} = 0,
 \ \Div \vc{v} = 0 \ \mbox{in}\ \D'((0,T) \times \Omega; R^N),
\ \mbox{for some}\ \tn{F} \in L^\infty((0,T) \times \Omega; R^{N \times N}_{{\rm sym},0} ),
\]
\[
\vc{v} \in C(Q; R^N),\ \tn{F}  \in C(Q; R^{N \times N}_{{\rm sym},0} ),
\]
\[
\sup_{(t,x) \in Q, t > \tau}
\left.
\frac{N}{2}\lambda_{\rm max} \left[ \frac{ (\vc{v} + \vc{h}[\vc{v}] ) \otimes (\vc{v} + \vc{h}[\vc{v}] ) }{r[\vc{v}]} - \tn{F} + \tn{M}[\vc{v}] \right]
- e[\vc{v}] < 0 \ \mbox{for any}\ 0 < \tau < T \vphantom{\frac{1}{2}}  \right\},
\]
where the symbol $\lambda_{\rm max}[\tn{A}]$ stands for the maximal eigenvalue of a symmetric matrix $\tn{A}$. We report the following result, see \cite[Theorem 2.1]{EF2015_1}:

\bProposition{E1}

Let the operators $\vc{h}$, $r$, $\tn{H}$, and $e$  be $Q-$continuous, where $Q \subset \left[ (0,T) \times \Omega \right]$ is an open set,
$
|Q| = |(0,T) \times \Omega|.
$
In addition, suppose that $r[\vc{v}] > 0$ and that the mapping $\vc{v} \mapsto 1/r[\vc{v}]$ is continuous in the sense specified in \eqref{E4a}.
Finally, assume that the set of subsolutions $X_0$ is \emph{non-empty} and \emph{bounded} in $L^\infty((0,T) \times \Omega; R^N)$.
Then problem \eqref{E1} -- \eqref{E3} admits infinitely many solutions.
\eP

In the remaining part of this section, we recast our problem (\ref{i1}), (\ref{i2}), (\ref{m4}) in the form \eqref{E1}--\eqref{E3} so that Proposition
\ref{PE1} may yield the conclusion claimed in Theorem \ref{Tm1}. To begin, we take
$
Q = (0,T) \times \Omega.
$

\subsection{Momentum decomposition}

We write the momentum $\vr \vu$ in the form
\[
\vr \vu = \vc{v} + \vc{V} + \Grad \Phi,
\]
where
\[
\Div \vc{v} = 0, \ \intO{ \Phi (t, \cdot) } = 0,\ \intO{ \vc{v}(t, \cdot) } = 0,\
\vc{V} = \vc{V}(t) \in R^2 - \mbox{a scalar function}.
\]
Next, we write the initial momentum by means of its Helmholtz decomposition as
$$
\vr_0 \vu_0 = \vc{v}_0 + \vc{V}_0 + \Grad \Phi_0, \ \Div \vc{v}_0 = 0, \ \intO{ \vc{v}_0 } = \intO{ \Phi_0 } = 0, \ \vc{V}_0 =
\frac{1}{|\Omega|} \intO{ \vr_0 \vu_0 }.
$$
Accordingly,
the continuity equation (\ref{i1}) reads
\bFormula{i3}
\partial_t \vr + \Del \Phi = 0 \ \mbox{in}\ (0,T) \times \Omega,
\eF
where the initial value of the density $\vr_0$ has been specified in (\ref{m4}).

Now, we can \emph{choose} $\vr = \vr(t,x) \in C^2([0,T] \times \Omega)$ in such a way that (\ref{i3}) holds for a certain potential
$\Phi$,
\[
\partial_t \vr (0, \cdot) = - \Delta \Phi_0,\ \ \Phi(0, \cdot) = \Phi_0,\
\intO{ \Phi(t, \cdot) } = 0 \ \mbox{for any}\ t \in [0,T].
\]
From now on, we will therefore assume that
\[
\vr \in C^2([0,T] \times \Omega), \ \Phi \in C^1([0,T]; C^3(\Omega; R^3))
\]
are fixed functions. Consequently, system (\ref{i1}), (\ref{i2}) reduces to
\begin{align}
\partial_t \vc{v} + \partial_t \vc{V} + \Div & \left[ \frac{ (\vc{v} + \vc{V} + \Grad \Phi ) \otimes (\vc{v} + \vc{V} + \Grad \Phi ) }{\vr}
+ \left( p(\vr)  + \partial_t \Phi \right) \tn{I} \right]
\nonumber\\
=& \left( \vc{v} + \vc{V} + \Grad \Phi \right)\left(1  - H\left( \frac{1}{\vr^2} \left| \vc{v} + \vc{V} + \Grad \Phi \right|^2 \right) \right) \label{i7}
\\
&-\vr \Grad K * \vr + \vr \psi * \left( \vc{v} + \vc{V} + \Grad \Phi \right) - \left( \vc{v} + \vc{V} + \Grad \Phi \right) \psi * \vr, \nonumber \\ \Div \vc{v} =&\,\, 0\,. \nonumber
\end{align}
Note that (\ref{i7}) still contains two unknowns, namely the
solenoidal field $\vc{v}$ and the scalar function $\vc{V}$.

\subsection{Kinetic energy}

Next,
we introduce the kinetic energy
\[
\frac{1}{2} \frac{ |\vc{v} + \vc{V} + \Grad \Phi |^2 }{\vr}
\]
and, motivated by (\ref{E2}), we look for solutions of (\ref{i7}) satisfying the following constraint:
\bFormula{i8}
\frac{1}{2} \frac{ |\vc{v} + \vc{V} + \Grad \Phi |^2 }{\vr} = e \equiv \Lambda - \frac{N}{2} \left( p(\vr) + \partial_t \Phi \right),
\eF
where $\Lambda = \Lambda(t)$ is a spatially homogeneous function to be determined below.
With this ansatz, relation (\ref{i7}) can be written as
\begin{align}
\partial_t \vc{v} + \partial_t \vc{V} + \Div \left[ \frac{ (\vc{v} + \vc{V} + \Grad \Phi ) \odot (\vc{v} + \vc{V} + \Grad \Phi ) }{\vr}
\right]
= &\left( \vc{v} + \vc{V} + \Grad \Phi \right)\left(1  - H \left( \frac{2}{\vr} e \right) \right)
\nonumber \\
&-\vr \Grad K * \vr + \vr \psi * \left( \vc{v} + \vc{V} + \Grad \Phi \right) \nonumber\\
&- \left( \vc{v} + \vc{V} + \Grad \Phi \right) \psi * \vr=: \Xi. \label{i9}
\end{align}

\subsection{The scalar function $\vc{V}$}

Our aim is to convert (\ref{i9}) to a ``balance law'' with a source term of zero mean. Accordingly, we fix $\vc{V}$ as the (unique) solution of the
following system of ordinary differential equations:
\begin{align}
\frac{{\rm d} \vc{V}}{{\rm d}t} + \vc{V} \frac{1}{|\Omega|} \left[ \intO{ \left( H \left( \frac{2}{\vr} e \right) - 1    \right)     }   \right]
= & - \frac{1}{|\Omega|} \intO{ \left( \vc{v} + \Grad \Phi \right) \left( H \left( \frac{2}{\vr} e \right) + \psi * \vr \right)  }\nonumber\\
&+ \frac{1}{|\Omega|} \intO{ \left( \vr \psi * \left( \vc{v} + \Grad \phi \right)  - \vr \Grad K * \vr  \right)},\label{i10}
\end{align}
with initial data $\vc{V}(0) = \vc{V}_0$, where we have used the identity
\[
\intO{   \vr } \intO{\psi} = \intO{ \psi * \vr }.
\]
Note that $\vc{V} = \vc{V}[\vc{v}]$ depends linearly on the solenoidal component $\vc{v}$ and also on the so far undetermined function $\Lambda$ through the quantity $e$.
Consequently, relation (\ref{i9}) reduces to
\begin{align*}
\partial_t \vc{v} + \Div \left[ \frac{ (\vc{v} + \vc{V} + \Grad \Phi ) \odot (\vc{v} + \vc{V} + \Grad \Phi ) }{\vr}
\right]
= \Xi - \frac{1}{|\Omega|} \intO{ \Xi }
\end{align*}
where the expression on the right-hand side has zero integral mean at any time $t$.
\subsection{Abstract problem}
Denoting
\[
\mathcal{G}[\vc{v}] = \Xi - \frac{1}{|\Omega|} \intO{ \Xi },
\]
we solve the elliptic problem
\bFormula{i12}
- \Div \left( \Grad \vc{w} + \Grad^t \vc{w} - \frac{2}{N} \Div \vc{w} \right) = \mathcal{ G} \ \mbox{in}\ \Omega \ \mbox{for any fixed}\ t \in [0,T]
\eF
and set
\[
\tn{M}[\vc{v}] = \Grad \vc{w} + \Grad^t \vc{w} - \frac{2}{N} \Div \vc{w}.
\]
Thus we have transformed the statement of Theorem \ref{Tm1} to the form compatible with (\ref{E1}--\ref{E3}), namely:

{\it
Find
$$
\vc{v} \in C_{\rm weak} ([0,T]; L^2(\Omega; R^3)), \ \Div \vc{v} = 0, \ \vc{v}(0, \cdot) = \vc{v}_0
$$
satisfying
$$
\partial_t \vc{v} + \Div \left[ \frac{ (\vc{v} + \vc{V}[\vc{v}] + \Grad \Phi ) \odot (\vc{v} + \vc{V}[\vc{v}] + \Grad \Phi ) }{\vr} + \tn{M}[\vc{v}]
\right] =0 ,
$$
together with the constraint
\bFormula{P3}
\frac{1}{2} \frac{ |\vc{v} + \vc{V}[\vc{v}] + \Grad \Phi |^2 }{\vr} = e \equiv \Lambda - \frac{N}{2} \left( p(\vr) + \partial_t \Phi \right).
\eF
}

\subsection{Proof of Theorem \ref{Tm1}}

It is a routine matter to check that the operators
\[
r = \vr, \ \vc{h} = \vc{V}[\vc{v}] + \Grad \Phi, \ \tn{M}[\vc{v}], \ \mbox{and}\ e
\]
comply with the hypotheses of Proposition \ref{PE1}, notably they are $Q$-continuous, with $Q = (0,T) \times \Omega$, see \cite[Section 2.3]{EF2015_1} for the relevant
discussion concerning $\tn{M}$.

Thus, in order to apply Proposition \ref{PE1}, it is enough to fix the function $\Lambda$ in (\ref{i8}) so that the set
of subsolutions $X_0$ is bounded and non-empty. To this end, we take
$\vc{v}_T = \vc{v}_0$ in (\ref{E3}) and consider
$\vc{v} = \vc{v}_0, \ \tn{F} = 0$
in (\ref{E6}). Consequently, it suffices to observe that $\Lambda = \Lambda(t)$ can be chosen in such a way that
\bFormula{EE1}
\sup_{(t,x) \in Q, t > \tau}
\frac{N}{2}\lambda_{\rm max} \left[ \frac{ (\vc{v_0} + \vc{V}[\vc{v_0}] + \Grad \Phi ) \otimes (\vc{v}_0 + \vc{V}[\vc{v_0}] + \Grad \Phi ) }{\vr} + \tn{M}[\vc{v}_0] \right]
- \Lambda(t) + \frac{N}{2} p(\vr) + \partial_t \Phi < 0
\eF
for any $0 < \tau < T$. As $\vc{v}_0$, $\vr$, $\Phi$ are fixed and $\vc{V}$, $\tn{M}$ are given by (\ref{i10}), (\ref{i12}) with
$H$ \emph{bounded}, there exists $\Lambda_0 > 0$ such that (\ref{EE1}) holds whenever
$$
\Lambda(t) \geq \Lambda_0 \ \mbox{for all}\ t \in [0,T]\,.
$$

\bRemark{R3}

Note that this is the only step in the proof of Theorem \ref{Tm1}, where we have used hypothesis (\ref{m3}).

\eR

Having fixed $\Lambda$ we can use the standard inequality (cf. \cite{DelSze3})
\[
\frac{1}{2} \frac{ |{\vc{h}} |^2 }{\tilde r} \leq
\frac{N}{2}\lambda_{\rm max} \left[ \frac{ {\vc{h}}  \otimes {\vc{h}} }{\tilde r} - \tn{H} \right]
\ \mbox{whenever}\ \tn{H} \in R^{N \times N}_{{\rm sym},0}
\]
to deduce that the set of subsolutions $X_0$ is bounded in $L^\infty((0,T) \times \Omega; R^3)$. Consequently, Proposition \ref{PE1} can be applied
to complete the proof of Theorem \ref{Tm1}.


\section{Dissipative solutions}
\label{d}

The solutions obtained via the method of convex integration delineated above suffer the essential difficulty that they typically do not obey the total energy balance associated to system (\ref{i1})-(\ref{i2}). Multiplying, formally, equation (\ref{i2}) on $\vu$ and integrating the resulting
expression over $\Omega$ we deduce
\begin{align}
\frac{{\rm d}}{{\rm d}t} \intO{ \left[ \frac{1}{2} \vr |\vu|^2 + P(\vr) \right] }
=& \intO{ \vr |\vu|^2 \left(1 - H\left(|\vu|^2 \right)\right) }-\intO{ \vr \vu \cdot \Grad K * \vr } \nonumber\\
&+ \intO{ \int_\Omega \psi (x - y) \vr(t,x) \vr(t,y) \vu(t,x)  \cdot \left( \vu(t,y) - \vu(t,x) \right) \ {\rm d}y }.\label{d1}
\end{align}
 If, in addition, we use the symmetry $K(Z) = K(-Z)$ and $\psi(x)=\psi(-x)$, we may write
\[
\intO{ \vr \vu \cdot \Grad K * \vr } = - \intO{ \Div (\vr \vu) K* \vr } = \frac{{\rm d}}{{\rm d}t} \frac{1}{2} \intO{ \vr K * \vr }
\]
and
\begin{align*}
\int_\Omega\int_\Omega \psi (x - y) \vr(t,x) \vr(t,y) \vu(t,x)  & \cdot \left( \vu(t,y) - \vu(t,x) \right) \ {\rm d}y \ {\rm d}x=\\
&-\intO{ \int_\Omega \psi (x - y) \vr(t,x) \vr(t,y) \left| \vu(t,y) - \vu(t,x) \right|^2 \ {\rm d}y};
\end{align*}
whence (\ref{d1}) reads
\begin{align}
\frac{{\rm d}}{{\rm d}t} \intO{ \left[ \frac{1}{2} \vr |\vu|^2 + P(\vr) + \frac{1}{2} \vr K* \vr  \right] }
= &\intO{ \vr |\vu|^2 \left(1 - H\left(|\vu|^2 \right)\right) }
\nonumber\\
&-\intO{ \int_\Omega \psi (x - y) \vr(t,x) \vr(t,y) \left| \vu(t,y) - \vu(t,x) \right|^2 \ {\rm d}y}\,.\label{d2}
\end{align}
In the framework of weak solutions, it is customary to impose (\ref{d2}) in the form of inequality
\begin{align}
\intO{ \left[ \frac{1}{2} \vr |\vu|^2 + P(\vr) + \frac{1}{2} \vr K* \vr  \right] (\tau, \cdot) }
\leq &
\intO{ \left[ \frac{1}{2} \vr_0 |\vu_0|^2 + P(\vr_0) + \frac{1}{2} \vr_0 K* \vr_0  \right] }
\label{d3}\\
&+ \int_0^\tau \intO{ \vr |\vu|^2 \left(1 - H\left(|\vu|^2 \right)\right) } \ \dt
\nonumber\\
&-\int_0^\tau\intO{ \int_\Omega \psi (x - y) \vr(t,x) \vr(t,y) \left| \vu(t,y) - \vu(t,x) \right|^2 \ {\rm d}y}\ \dt\,
\nonumber
\end{align}
for a.e. $\tau \in (0,T)$, as an admissibility criterion. The weak solutions satisfying (\ref{d3}) are called \emph{dissipative weak solutions}. Our ultimate goal is to identify a class of initial data for which problem (\ref{i1}), (\ref{i2}), (\ref{m4}) admits infinitely many \emph{dissipative} weak solutions. Revisiting the proof of Theorem \ref{Tm1} in Section \ref{r} we take
\[
\vr (t,\cdot) = \vr_0 \in C^2(\Omega)\ \mbox{for all}\ t \in [0,T],
\]
with the associated potential $\Phi = 0$ and the initial momentum
\[
\vr_0 \vu_0 = \vc{v}_0, \ \vc{v}_0 \in C^3(\Omega), \ \Div \vc{v}_0 = 0, \ \intO{ \vc{v}_0 } = 0,
\]
and, accordingly, $\vc{V}_0 = 0$. Going back to (\ref{d2}) and using (\ref{P3}), we compute
\[
\frac{{\rm d}}{{\rm d}t} \intO{ \left[ \frac{1}{2} \vr |\vu|^2 + P(\vr) + \frac{1}{2} \vr K* \vr  \right] } =
|\Omega| \frac{{\rm d}}{{\rm d}t} \Lambda,
\]
while the right-hand side of (\ref{d2}) can be trivially estimated by
\[
\intO{ \vr |\vu|^2 } \leq c( 1 + \Lambda)\,,
\]
according to the assumptions on $H$. Consequently, taking
$
\Lambda(t) = \Lambda_0 + \exp(- \lambda t) \ \mbox{with}\ \lambda > 0 \ \mbox{large enough},
$
the weak solutions resulting from the procedure elaborated in Section \ref{r} will satisfy
\begin{align}\label{d4}
\frac{{\rm d}}{{\rm d}t} \int_\Omega \left[ \frac{1}{2} \vr |\vu|^2 \right. &\left.+ P(\vr) + \frac{1}{2} \vr K* \vr  \right]  \ {\rm d}x
\nonumber\\
&\leq \intO{ \vr |\vu|^2 \left(1 - H\left(|\vu|^2 \right)\right) }-\intO{ \int_\Omega \psi (x - y) \vr(t,x) \vr(t,y) \left| \vu(t,y) - \vu(t,x) \right|^2 \ {\rm d}y}
\end{align}
In case $\psi$ is not symmetric, one can always keep the corresponding term in the right-hand side of \eqref{d3} for straightforward estimates due to the smoothness of $\psi$. We leave the details to the reader. Finally, to deduce (\ref{d3}) from (\ref{d4}), it remains to find a family of subsolutions for which the energy does not ``jump up'' at the initial time $\tau = 0$. The necessary piece of information is provided by the following result proved in \cite[Theorem 6.1]{EF2015_1}:

\bProposition{r1}

In addition to the hypotheses of Proposition \ref{PE1}, suppose that
\bFormula{HYP}
\left| \left\{ x \in \Omega \ \Big| (t,x) \in Q \right\} \right| = |\Omega| \ \mbox{for any}\ 0 < t < T.
\eF

Then there exists a set of times $\mathcal{R} \subset (0,T)$ dense in $(0,T)$ such that for
any $\tau \in \mathcal{R}$ there is $\vc{v} \in \Ov{X}_0$ with the following properties:
\begin{itemize}
\item[i)]
$\vc{v} \in C_b(\left[ (0,\tau) \cup (\tau, T) \times {\Omega} \right] \cap Q; R^N) \cap C_{\rm weak}([0,T]; L^2(\Omega; R^N)) ,\
\vc{v}(0, \cdot) = \vu_0 , \ \vc{v}(T, \cdot) = \vu_T$.

\item[ii)]
$\partial_t \vc{v} + \Div \tn{F} = 0,
 \ \Div \vc{v} = 0 \ \mbox{in}\ \D'((0,T) \times \Omega; R^N)$,
for some $\tn{F} \in C_b( \left[ (0,\tau) \cup (\tau, T) \times \Omega \right] \cap Q; R^{3 \times 3}_{{\rm sym},0} )$.
\item[iii)]
$$\frac{N}{2}\lambda_{\rm max} \left[ \frac{ (\vc{v} + h[\vc{v}] ) \otimes (\vc{v} + h[\vc{v}] ) }{r[\vc{v}]} - \tn{F} + \tn{M}[\vc{v}] \right]
 < e[\vc{v}] \vphantom{\frac{1}{2}} \ \mbox{in}\ \left[ (0,\tau)   \times \Omega \right] \cap Q\,.$$
\item[iv)]
$$\sup_{(t,x) \in Q, t > \tau + s}
\frac{N}{2}\lambda_{\rm max} \left[ \frac{ (\vc{v} + h[\vc{v}] ) \otimes (\vc{v} + h[\vc{v}] ) }{r[\vc{v}]} - \tn{F} + \tn{M}[\vc{v}] \right]
 - e[\vc{v}] < 0 \ \mbox{for any}\ 0 < s < T- \tau\,.$$
\item[v)]
$$\frac{1}{2} \intO{ \frac{ |\vc{v} + \vc{h}[\vc{v}] |^2 }{r[\vc{v}]} ({\tau}, \cdot) } = \intO{ e[\vc{v}]({\tau}, \cdot) }.$$

\end{itemize}
\eP

Since in our setting $Q = (0,T) \times \Omega$ hypothesis (\ref{HYP}) is obviously satisfied, we may combine the final arguments of Theorem \ref{Tm1} with
(\ref{d4}) and Proposition \ref{Pr1} to conclude the proof of Theorem \ref{Td1}.

\bRemark{r3}
Observe that $\vr = \vr_0$ has been chosen constant in time in this section so that the suitable initial data provided by Proposition
\ref{Pr1} correspond to the solutions of \eqref{i1}, \eqref{i2} and \eqref{m4}, with
\[
\vr(0, \cdot) = \vr_0.
\]
\eR


\section{Relative energy and weak-strong uniqueness}
\label{WSU}

Our ultimate goal is to show that a dissipative and a strong solution originating from the same initial data coincide as long as the latter exists. To this end, we revoke the method proposed by Dafermos \cite{Daf4} and later elaborated in \cite{FeJiNo}, based on the concept of \emph{relative energy}.  We introduce the relative energy functional
$$
\mathcal{E} \left( \varrho , \vu \ \Big|\ r, \vc{U} \right) = \intO{ \left[ \frac{1}{2} \varrho |\vu - \vc{U}|^2  +P(\vr)-P'(r)(\vr-r)-P(r)\right] },
$$
where $[\varrho, \vu]$ is a dissipative weak solution  and $[r, \vc{U}]$ are sufficiently smooth functions in
$[0,T] \times \Omega$, $r> 0$.
Evoking the arguments of \cite[Section 3]{FeJiNo} we derive the relative energy inequality
\begin{equation} \label{rei}
\begin{split}
\left[ \mathcal{E} \left( \varrho , \vu \ \Big|\ r, \vc{U} \right) \right]_{t = 0}^{t = \tau} & + \frac{1}{2} \left[ \intO{ (r-\varrho) (K*(r-\vr))}  \right]_{t = 0}^{t = \tau}
\\
&+\frac12 \int_0^\tau \intO{ \int_{\Omega} \psi (x - y) \vr(t,x) \vr(t,y) \Big| (\vu(t,y) - \vc{U}(t,y))-(\vu(t,x)-\vc{U}(t,x) )\Big|^2 \ {\rm d}y } \ \dt \\ & \leq \int_0^\tau \intO{ \mathcal{R} \left( \vr, \vu \Big| r, \vc{U} \right) } \ \dt
\end{split}
\end{equation}
with the remainder term
\begin{equation} \label{reirem}
\begin{split}
{\mathcal R}\left(\vr, \vu \ \Big| \ r,\vc{U} \right)=& -\intO{\vr(\partial_t\vc{U}+\vu \cdot \Grad \vc{U}) \cdot (\vu-\vc{U})}\\
&-\intO{ \vr (\Grad K * r) \cdot \Big( \vu - \vc{U} \Big) } + \intO{ (\vr - r) \left( \Grad K* (\vr - r) \right) \cdot \vc{U} }  \\
&+\intO{(r-\vr)\partial_tP'(r)+\Grad P'(r) \cdot (r\vc{U}-\vr\vu)-\Div\vc{U}(p(\vr)-p(r))}\\
&+\intO{\left( 1 - H\left( |\vu |^2 \right) \right) \vr \vu \cdot (\vu-\vc{U})}\\
&+ \intO{ \int_{\Omega} \psi (x-y) \vr(t,x)(\vu(t,x) - \vc{U}(t,x) ) \cdot ( \vc{U}(t,y) - \vc{U}(t,x) ) ( \vr(t,y) - r(t,y) ) \ {\rm d}y } \\
&- \intO{ \int_{\Omega} \vr(t,x) \psi(x- y) ( \vc{U}(t,y) - \vc{U}(t,x) ) r(t,y) \ {\rm d}y \cdot (\vc{U}(t,x) - \vc{u}(t,x) ) \ } \ \dt.
\end{split}
\end{equation}
Relations (\ref{rei}), (\ref{reirem}) are valid for any dissipative weak solution $[\vr, \vu]$ and any pair of smooth ``test'' functions $[r, \vc{U}]$, $r > 0$, cf. \cite[Section 3]{FeJiNo}.
With (\ref{rei}), (\ref{reirem}) at hand, we are ready to prove Theorem \ref{Te1}.

Similarly to \cite{FeJiNo}, the proof consists in taking $[r, \vc{U}]$ as test functions in (\ref{rei}) and applying a Gronwall type argument. Since $r$ is smooth and bounded below, we may find two
positive constant $\underline{\vr}$, $\Ov{\vr}$ such that
\bFormula{rbounds}
0 < \underline{\vr} < \min_{[0,T] \times \Omega} r \leq \max_{[0,T] \times \Omega} r < \Ov{\vr} < \infty.
\eF
Now it is convenient to introduce the essential and residual component of a measurable function $h$ as
\[
h = h_{\rm ess} + h_{\rm res},\ h_{\rm ess} = \chi (\vr) h,\ \chi \in \DC(0, \infty), \ 0 \leq \chi \leq 1,\ \chi |_{[\underline{\vr}, \Ov{\vr}]} = 1.
\]

In view of hypothesis (\ref{tlak}) and (\ref{rbounds}), it can be checked that
\begin{equation} \label{coerc}
\begin{split}
\left[ \tfrac{1}{2} \varrho |\vu - \vc{U}|^2\right. & \left.+P(\vr)-P'(r)(\vr-r)-P(r)\right]  \\
&\geq c(\underline{\vr}, \Ov{\vr}) \left( \vr | \vu - \vc{U} |^2 + \left| \left[ \vu - \vc{U} \right]_{\rm ess} \right|^2 +  \left| \left[ \vr - r \right]_{\rm ess} \right|^2
+ \left( 1_{\rm res} + [p(\vr)]_{\rm res} + [\vr \log^+(\vr)]_{\rm res} \right) \right), \ c > 0,
\end{split}
\end{equation}
cf. \cite[formula (4.1)]{FeJiNo}. From now on, $c$ will denote generic constants that will change from line to line and that we avoid to specify fully. A direct consequence of this property and the fact that the total mass of the weak solution, 
\[
M = \intO{ \vr(t,\cdot) } = \intO{ \vr_0 }
\]
is a constant of motion, is the following estimate:
\begin{equation}\label{kkk}
\left\| \left[ \vr - r \right]_{\rm ess} \right\|_{L^1(\Omega)} + \left\| \left[ \vr - r \right]_{\rm res} \right\|_{L^1(\Omega)}+\left\| \left[ \vr - r \right]_{\rm ess} \right\|_{L^2(\Omega)}
\leq c(\underline{\vr}, \Ov{\vr} , M)  \, \mathcal{E} \left( \varrho , \vu \ \Big|\ r, \vc{U} \right)^{\frac{1}{2}}\,.
\end{equation}
Indeed the estimate of the essential component follows directly from (\ref{coerc}) and the embedding $L^2(\Omega) \hookrightarrow L^1(\Omega)$, while, 
using (\ref{coerc}) once more, we have
\[
\left\| [\vr - r]_{\rm res} \right\|_{L^1(\Omega)} \leq  (M + |\Omega| \Ov{\vr})^{1/2} \left\| [\vr - r]_{\rm res} \right\|^{1/2}_{L^1(\Omega)} 
\leq c(\underline{\vr}, \Ov{\vr} , M)  \, \mathcal{E} \left( \varrho , \vu \ \Big|\ r, \vc{U} \right)^{\frac{1}{2}}
\]
as desired.

The proof of Theorem \ref{Te1} will be carried over in several steps.
\color{black}

\medskip

{\bf Step 1.} Plugging the strong solution $[r, \vc{U}]$ in (\ref{reirem}) we write
\[
\vr \left( \partial_t \vc{U} + \vu \cdot \Grad \vc{U} \right) \cdot (\vu - \vc{U} ) = \vr (\vu - \vc{U} ) \cdot \Grad \vc{U} \cdot (\vc{u} - \vc{U}) +
\vr \left( \partial_t \vc{U} + \vc{U} \cdot \Grad \vc{U} \right),
\]
where the first term on the right-hand side can be ``absorbed'' in \eqref{rei} by a direct $L^2$-estimate while the second can be rewritten by means of (\ref{i2}) satisfied by $\vc{U}$. Accordingly,  the relative energy inequality (\ref{rei}) gives rise to
\begin{equation} \label{rei1}
\begin{split}
&\left[ \mathcal{E} \left( \varrho , \vu \ \Big|\ r, \vc{U} \right) \right]_{t = 0}^{t = \tau}  +
\frac{1}{2} \left[ \intO{ (r-\varrho) (K*(r-\vr))}  \right]_{t = 0}^{t = \tau}
\\
& +\frac12  \int_0^\tau \intO{ \int_{\Omega} \psi (x - y) \vr(t,x) \vr(t,y) \Big| (\vu(t,y) - \vc{U}(t,y))-(\vu(t,x)-\vc{U}(t,x) )\Big|^2 \ {\rm d}y } \ \dt \\
& \leq c \int_0^\tau \mathcal{E} \left( \varrho , \vu \ \Big|\ r, \vc{U} \right) \ \dt +\int_0^\tau
 \intO{ (\vr - r) \left( \Grad K* (\vr - r) \right) \cdot \vc{U}  } \ \dt \\
& \quad+ \int_0^\tau \intO{ \left[ (r-\vr)\partial_tP'(r)+\Grad P'(r) \cdot (r\vc{U}-\vr\vu)-\Div\vc{U}(p(\vr)-p(r)) + \frac{\vr}{r} \Grad p(r) \cdot (\vu - \vc{U})   \right] } \ \dt \\
&\quad+ \int_0^\tau \intO{\left( 1 - H\left( |\vu |^2 \right) \right) \vr \vu \cdot (\vu-\vc{U})} \ \dt - \int_0^\tau \intO{\left( 1 - H\left( |\vc{U} |^2 \right) \right) \vr \vc{U} \cdot (\vu-\vc{U})} \ \dt\\
&\quad+ \int_0^\tau \intO{ \int_{\Omega} \psi (x-y) \vr(t,x)(\vu(t,x) - \vc{U}(t,x) ) \cdot ( \vc{U}(t,y) - \vc{U}(t,x) ) ( \vr(t,y) - r(t,y) ) \ {\rm d}y } \ \dt.
\end{split}
\end{equation}

{\bf Step 2.} Next, in view of hypothesis (\ref{HH1}), we may write
$
H(Z) = H_1(Z) + H_2(Z),
$
where $H_1$ has compact support and $H_2$ is a non-decreasing function. We have
\[
\begin{split}
\left|\int_0^\tau \intO{\left( 1 - H_1\left( |\vu |^2 \right) \right) \vr \vu \cdot (\vu-\vc{U})} \ \dt \right. &-\left. \int_0^\tau \intO{\left( 1 - H_1\left( |\vc{U} |^2 \right) \right) \vr \vc{U} \cdot (\vu-\vc{U})} \ \dt \right| \\
& \leq c \int_0^\tau \intO{ \vr | \vu - \vc{U} |^2 } \ \dt \leq c \int_0^\tau \mathcal{E} \left( \varrho , \vu \ \Big|\ r, \vc{U} \right) \ \dt,
\end{split}
\]
while
\[
\begin{split}
\int_0^\tau \intO{&\left( 1 - H_2\left( |\vu |^2 \right) \right) \vr \vu \cdot (\vu-\vc{U})} \ \dt - \int_0^\tau \intO{\left( 1 - H_2\left( |\vc{U} |^2 \right) \right) \vr \vc{U} \cdot (\vu-\vc{U})} \ \dt \\
& = \int_0^\tau \intO{ \vr |\vu - \vc{U}|^2 } \ \dt - \int_0^\tau \intO{ \left( H_2\left( |\vc{u} |^2 \right) \vu - H_2 \left( |\vc{U} |^2 \right) \vc{U} \right)
\cdot (\vu - \vc{U}) } \ \dt\\
&\leq \int_0^\tau \intO{ \vr |\vu - \vc{U}|^2 } \ \dt \leq 2 \int_0^\tau \mathcal{E} \left( \varrho , \vu \ \Big|\ r, \vc{U} \right) \ \dt.
\end{split}
\]
Next, by virtue of (\ref{coerc}), we deduce
\[
\begin{split}
\left| \int_0^\tau \intO{ \right.&\left.\left[ (r-\vr)\partial_tP'(r)+\Grad P'(r) \cdot (r\vc{U}-\vr\vu)-\Div\vc{U}(p(\vr)-p(r)) + \frac{\vr}{r} \Grad p(r) \cdot (\vu - \vc{U})   \right] } \ \dt \right| \\
& =  \left| \int_0^\tau \intO{ \Div \vc{U} \Big( p(\vr) - p'(r)(\vr - r) - p(r) \Big) } \ \dt \right|\leq c \int_0^\tau \mathcal{E} \left( \varrho , \vu \ \Big|\ r, \vc{U} \right) \ \dt;
\end{split}
\]
whence (\ref{rei1}) reduces to
\begin{equation} \label{rei2}
\begin{split}
\left[ \mathcal{E} \right.&\left.\left( \varrho , \vu \ \Big|\ r, \vc{U} \right) \right]_{t = 0}^{t = \tau} +
\frac{1}{2} \left[ \intO{ (r-\varrho) (K*(r-\vr))}  \right]_{t = 0}^{t = \tau}\\
&+ \frac12 \int_0^\tau \intO{ \int_{\Omega} \psi (x - y) \vr(t,x) \vr(t,y) \Big| (\vu(t,y) - \vc{U}(t,y))-(\vu(t,x)-\vc{U}(t,x) )\Big|^2 \ {\rm d}y } \ \dt \\
\leq &\, c \int_0^\tau \mathcal{E} \left( \varrho , \vu \ \Big|\ r, \vc{U} \right) \ \dt + \int_0^\tau \intO{ (\vr - r) \left( \Grad K* (\vr - r) \right) \cdot \vc{U}  } \ \dt \\
&+ \int_0^\tau \intO{ \int_{\Omega} \psi (x-y) \vr(t,x)(\vu(t,x) - \vc{U}(t,x) ) \cdot ( \vc{U}(t,y) - \vc{U}(t,x) ) ( \vr(t,y) - r(t,y) ) \ {\rm d}y } \ \dt.
\end{split}
\end{equation}

\medskip

{\bf Step 3.} Using H\" older's and Young's inequalities, we deduce
\begin{equation*}
\begin{split}
\Big{|}  \intO{ &\int_{\Omega} \psi (x-y) \vr(t,x)(\vu(t,x) - \vc{U}(t,x) ) \cdot ( \vc{U}(t,y) - \vc{U}(t,x) ) ( \vr(t,y) - r(t,y) ) \ {\rm d}y } \Big{|}\\
&\le c\|\vc{U}\|_{L^\infty(\Omega)}\left(\intO{\vr |\vu -\vc{U}|^2}\right)^{\frac{1}{2}}\ \left(\intO{ \vr(t,x)\left(\int_\Omega \psi(x-y)(\vr(t,y)-r(t,y)) \ {\rm d}y}\right)^2 \ \dx \right)^\frac{1}{2} \\
&\le c\, \mathcal{E} \left( \varrho , \vu \ \Big|\ r, \vc{U} \right)  +c \intO{\vr(t,x)} \left(\sup_{x\in\Omega}\int_\Omega \psi(x-y)(\vr(t,y)-r(t,y))\ {\rm d}y\right)^2 \\
& \leq c\, \mathcal{E} \left( \varrho , \vu \ \Big|\ r, \vc{U} \right) + c \left(\sup_{x\in\Omega}\int_\Omega \psi(x-y)(\vr(t,y)-r(t,y))\ {\rm d}y\right)^2.
\end{split}
\end{equation*}
Now, we can use \eqref{kkk} to get
\[
\begin{split}
\left| \int_\Omega \psi(x-y)(\vr(t,y)-r(t,y))\ {\rm d}y \right| \leq c \left\| \vr - r \right\|_{L^1(\Omega)}
\leq c \,\mathcal{E} \left( \varrho , \vu \ \Big|\ r, \vc{U} \right)^{\frac{1}{2}}\,.
\end{split}
\]

In order to deal with the potential term, we simply use \eqref{kkk} to obtain
\begin{align*}
\left| \int_0^\tau
\intO{ (\vr - r) \left( \Grad K* (\vr - r) \right) \cdot \vc{U}  } \ \dt \right|
&\leq c \int_0^\tau \left\| \Grad K* (\vr - r) \right\|_{L^\infty(\Omega)} \left\| \vr - r\right\|_{L^1(\Omega)} \dt \\
&\leq c \int_0^\tau \left\| \vr - r \right\|_{L^1(\Omega)}^2 \dt \leq c \int_0^\tau \mathcal{E} \left( \varrho , \vu \ \Big|\ r, \vc{U} \right) \ \dt.
\end{align*}
\color{black}

Summing up, we reduced (\ref{rei2}) to
\begin{equation} \label{rei3}
\begin{split}
\left[ \mathcal{E} \left( \varrho , \vu \ \Big|\ r, \vc{U} \right) \right]_{t = 0}^{t = \tau} +
\frac{1}{2} \left[ \intO{ (r-\varrho) (K*(r-\vr))}  \right]_{t = 0}^{t = \tau}
\leq c \int_0^\tau \mathcal{E} \left( \varrho , \vu \ \Big|\ r, \vc{U} \right) \ \dt.
\end{split}
\end{equation}

\medskip

{\bf Step 4.} Finally, we write
\[
\begin{split}
\frac{1}{2} \left[ \intO{ (r-\varrho) (K*(r-\vr))}  \right]_{t = 0}^{t = \tau} &= \int_0^\tau \intO{ (r-\varrho) (K*(\partial_t r-\partial_t \vr))}  \ \dt\\
& = \int_0^\tau \intO{ (\varrho - r ) (K*( \Div(r \vc{U}) -\Div (\vr \vu) ))}  \ \dt \\
& = \int_0^\tau \intO{ (r \vc{U} - \vr \vu ) \cdot \Grad K* (r - \vr) } \ \dt.
\end{split}
\]
Of course, this is a bit formal as $\vr$, $\vu$ represent only a weak solution but since the kernel $K$ is smooth the procedure can be justified by a density argument.
Furthermore, we have
\[
\int_0^\tau \intO{ (r \vc{U} - \vr \vu ) \cdot \Grad K* (r - \vr) } \ \dt =
\int_0^\tau \intO{ (r  - \vr ) \vc{U} \cdot \Grad K* (r - \vr) } \ \dt + \int_0^\tau \intO{ \vr (\vc{U} - \vu) \cdot \Grad K* (r - \vr) } \ \dt ,
\]
where the former integral on the right-hand side can be handled exactly as in {\bf Step 4}, while the latter reads
\[
\left| \int_0^\tau \!\!\!\intO{ \vr (\vc{U} - \vu) \cdot \Grad K* (r - \vr) } \ \dt \right| \leq c \left( \int_0^\tau \!\!\!\intO{ \vr |\vc{U} - \vu|^2 } \ \dt
+ \int_0^\tau \!\!\left[ \intO{ \vr } \left\| K* (r - \vr) \right\|^2_{L^\infty(\Omega)} \right] \ \dt \right),
\]
where the last term can be handled as in {\bf Step 4}.

Going back to (\ref{rei3}) we may therefore use Gronwall's lemma to obtain the desired conclusion.
\color{black}


\section{Pressureless system}\label{Poi}

In the proof of Theorem \ref{Te1}, the presence of the pressure $p$ plays a crucial role when dominating the convolution terms. On the other hand, certain models
do not count with $p$. In this final part, we show that the repulsive effect of the pressure can be substituted by a weak interaction term represented
by the Poisson kernel. Specifically, we consider the following system
\begin{align}
\partial_t \vr(t,x) + \Div \Big( \vr(t,x) \vu(t,x) \Big) =& \, 0,
\label{Pi1}\\
\partial_t \Big( \vr(t,x) \vu(t,x) \Big) + \Div \Big( \vr(t,x) \vu(t,x)  \otimes \vu(t,x) \Big) = & \left( 1 - H\left( |\vu(t,x) |^2 \right) \right) \vr(t,x) \vu(t,x) - \vr(t,x) \Grad \Phi_\varrho(t,x)
\nonumber\\
&- \vr(t,x) \int_{\Omega} \Grad K(x - y) \vr(t,y) \ {\rm d}y
\label{Pi2}\\
&+
\vr(t,x) \int_{\Omega} \psi (x-y) \Big( \vu(t,y) - \vu(t,x) \Big) \vr(t,y) \ {\rm d}y ,
\nonumber \\
-\Delta_x\Phi_\vr(t,x)=&\vr(t,x)-\overline\vr,
\label{Pi3}
\end{align}
where $\Ov{\vr} = \intO{ \vr }$ represent the total mass. Note that the solution to \eqref{Pi3} can be represented by means of a convolution
$$\Phi_\vr(t,x)=(\mathcal{P} *\vr)(t,x),$$
where $\mathcal{P}$ is the Poisson kernel.

 For strong solutions $[r, \vc{U}]$ it is convenient to rewrite  system \eqref{Pi1}--\eqref{Pi3}
as follows
\begin{align}
\partial_t r(t,x) + \Div \Big( r(t,x) \vc{U}(t,x) \Big) =& \, 0,
\label{i1-Pr}\\
\partial_t  \vc{U}(t,x)  + \vc{U}(t,x)\cdot\Grad \vc{U}(t,x) = & \left( 1 - H\left( |\vc{U}(t,x) |^2 \right) \right)  \vc{U}(t,x)
- \Grad\Phi_r(t,x)
\label{i2-Pr}\\
&+
 \int_{\Omega} \psi (x-y) \Big( \vc{U}(t,y) - \vc{U}(t,x) \Big) r(t,y) \ {\rm d}y ,
\nonumber\\
-\Delta_x\Phi_r(t,x)=&r(t,x)-\overline r.
\label{i3-Pr}
\end{align}
In the case there appears no vacuum, the systems \eqref{Pi1}--\eqref{Pi3}  and \eqref{i1-Pr}-\eqref{i3-Pr}
are equivalent for smooth solutions. However, choosing $[r,\vc{U}]$ as a solution to the second systems does not require to exclude compactly supported smooth solutions. The issue of relative entropy inequality for Euler-Poisson system was considered in~\cite{tzavaras}. In the present consideration we include also the alignment term. 
We observe that the repulsive force imposed on the system by the Poisson kernel can substitute the effect of the pressure in Theorem \ref{Te1}.

\bTheorem{e2}
Let $T > 0$ be given and let $N=2,3$. Suppose that  $\psi$, $K$ and $H$ satisfy the hypotheses \eqref{m1}--\eqref{m3}.  Let $[\varrho,\vu]$ be a dissipative weak solution to system \eqref{Pi1}--\eqref{Pi3} in $(0,T) \times \Omega$ with initial data $\varrho_0, \vu_0$ and $\vr_0\ge0$
 and let  $[r, \vc{U}]$ be a Lipschitz (strong) solution  to system \eqref{i1-Pr}-\eqref{i3-Pr} with\[
\varrho_0 = r(0, \cdot),\
\vu_0 = \vc{U}(0, \cdot).
\]
Then
\[
\varrho = r,\ \varrho\vu = r\vc{U} \ \mbox{a.e. in}\ (0,T) \times \Omega.
\]

\eT

{\bf Proof:} The proof mimicks the steps undertaken in the proof of Theorem \ref{Te1}. In the present setting, however, it is convenient to
introduce a relative energy in the form
\bFormula{e3-P}
\mathcal{E} \left( \varrho , \vu \ \Big|\ r, \vc{U} \right) = \intO{ \left[ \frac{1}{2} \varrho |\vu - \vc{U}|^2  +\frac{1}{2}(r-\varrho) (\mathcal{P}*(r-\vr))\right] },
\eF
where the new term represented by the Poisson kernel is equivalent to the negative Sobolev norm
\bFormula{norm}
\intO{ \frac{1}{2}(r-\varrho) (\mathcal{P}*(r-\vr)) } = \| \vr - r \|^2_{W^{-1,2}(\Omega)}.
\eF

To begin, note that {\bf Steps 1, 2} in the proof of Theorem \ref{Te1} pass without any changes in the present setting.

The first part of {\bf Step 3} can be carried over as soon as we realize that
\begin{equation} \label{regu}
\left\| \psi * (\vr - r) \right\|_{L^\infty(\Omega)} \leq c \left\| \vr - r \right\|_{W^{-1,2}(\Omega)}
\end{equation}
as the kernel $\psi$ is smooth.

In the second part of {\bf Step 3} and in {\bf Step 4}, we may use an inequality similar to (\ref{regu}) to handle the terms containing convolution with the kernel $K$. Thus the only
new ingredient is to control the term
\[
\intO{ (\vr - r) \left( \Grad \mathcal{P} * (\vr - r) \right) \cdot \vc{U}  } \ \dt.
\]
To this end, we proceed as in \cite{tzavaras} observing that
\[
(\vr - r) \left( \Grad \mathcal{P} * (\vr - r) \right) \cdot \vc{U}  =
\frac{1}{2} \Grad | \Grad \Phi_r - \Grad \Phi_\vr |^2 - \Div \left( \Grad (\Phi_r - \Phi_\vr) \otimes \Grad (\Phi_r - \Phi_\vr) \right),
\]
therefore
\begin{equation*}
\begin{split}
\int_0^\tau\intO{&\vc{U}(\vr-r)(\Grad \mathcal{P})*(\vr-r)}\\
&=-\int_0^\tau\intO{\left(\Div\vc{U} \left(\frac{1}{2}|\Grad(\Phi_r-\Phi_\vr)|^2 \right)-\Grad\vc{U} \cdot \left(\Grad(\Phi_r-\Phi_\vr)\otimes\Grad(\Phi_r-\Phi_\vr)\right)\right)}\ \dt.
\end{split}
\end{equation*}
By virtue of (\ref{e3-P}), (\ref{norm}), we may infer that
\[
\left| \intO{ (\vr - r) \left( \Grad \mathcal{P} * (\vr - r) \right) \cdot \vc{U}  } \ \dt \right| \leq c \,\mathcal{E} \left( \varrho , \vu \ \Big|\ r, \vc{U} \right),
\]
which completes the proof exactly as in Section \ref{WSU}. 
\qed

\begin{Remark}
In the contrast to Theorem~\ref{Te1} the velocities $\vu$ and $\vc{U}$ coincide up to the vacuum. 
\end{Remark}

\def\cprime{$'$} \def\ocirc#1{\ifmmode\setbox0=\hbox{$#1$}\dimen0=\ht0
  \advance\dimen0 by1pt\rlap{\hbox to\wd0{\hss\raise\dimen0
  \hbox{\hskip.2em$\scriptscriptstyle\circ$}\hss}}#1\else {\accent"17 #1}\fi}


\section*{Acknowledgments}
J.A.C. was partially supported by the project MTM2011-27739-C04-02 DGI (Spain), from the Royal Society by a Wolfson Research Merit Award and by the EPSRC grant EP/K008404/1. The research of E.F. leading to these results has received funding from the European Research Council under the European Union's Seventh Framework Programme (FP7/2007-2013)/ ERC Grant Agreement 320078 (The Institute of Mathematics of the Academy of Sciences of the Czech
Republic is supported by RVO:67985840).  The research of A.\'S.-G. has received funding from the  National Science Centre, DEC-2012/05/E/ST1/02218. The research of P. G. has received funding from the  National Science Centre, Poland, 2014/13/B/ST1/03094.

%
%
%
%

\end{document}